\documentclass[11pt]{article}

\def\mathbb{\Bbb}
\usepackage{amsfonts,latexsym,amstext}
\usepackage{amsmath}
\usepackage[active]{srcltx}
\usepackage{amssymb}

\newtheorem{theorem}{Theorem}[section]
\newtheorem{lemma}[theorem]{Lemma}
\newtheorem{proposition}[theorem]{Proposition}
\newtheorem{definition}{Definition}[section]

\newtheorem{hypothesis}[theorem]{Hypothesis}

\newtheorem{remark}[theorem]{Remark}
\newtheorem{corollary}[theorem]{Corollary}
\numberwithin{equation}{section}
\def\qed{{\hfill\hbox{\enspace${ \square}$}} \smallskip}
\def\sqr#1#2{{\vcenter{\vbox{\hrule height .#2pt \hbox{\vrule
 width .#2pt height#1pt \kern#1pt \vrule
width .#2pt} \hrule height .#2pt}}}}
\def\square{\mathchoice\sqr54\sqr54\sqr{4.1}3\sqr{3.5}3}

\def\ds{\begin{displaystyle}}
\def\eds{\end{displaystyle}}
\def\dis{\displaystyle }
\def\<{\langle }
\def\>{\rangle }

\def\dim{\noindent \hbox{{\bf Proof.} }}
\def\R{\mathbb R}
\def\N{\mathbb N}

\def\E{\mathbb E}
\def\P{\mathbb P}

\def\calg{{\cal G}}

\def\call{{\cal L}}

\def\calo{{\cal O}}

\begin{document}

\title{HJB equations in infinite dimensions under weak regularizing properties}
\date{}
 \author{
 Federica Masiero\\
 Dipartimento di Matematica e Applicazioni, Universit\`a di Milano Bicocca\\
 via Cozzi 55, 20125 Milano, Italy\\
 e-mail: federica.masiero@unimib.it
 }

\maketitle  
\begin{abstract}
We solve in mild sense Hamilton Jacobi Bellman equations, both in an infinite dimensional Hilbert space
and in a Banach space,
with lipschitz Hamiltonian and lipschitz continuous
final condition, and asking only a weak regularizing property on the transition semigroup of the corresponding state equation.
The results are applied to solve stochastic optimal control problems; the models we can treat include a controlled stochastic heat equation in space dimension one and with control and noise on a subdomain.
\end{abstract}

\section{Introduction}

In this paper we study semilinear Kolmogorov equations
in an infinite dimensional Hilbert space $H$, as well as in a
Banach space $E$, in particular Hamilton Jacobi Bellman equations.
More precisely, let us consider the following equation
\begin{equation}
\left\{
\begin{array}
[c]{l}%
-\frac{\partial v}{\partial t}(t,x)=\mathcal{L}_tv\left(  t,x\right)
+\psi\left(  \nabla v\left(  t,x\right)B  \right)+l(t,x)  ,\text{ \ \ \ \ }t\in\left[  0,T\right]
,\text{ }x\in H\text{ or }x\in E\\
v(T,x)=\phi\left(  x\right).
\end{array}
\right.  \label{Kolmo intro}%
\end{equation}
The second order differential operator $\call_t $ is the generator of the transition
semigroup $P_t$ related to the following perturbed Ornstein-Uhlenbeck process
\begin{equation}
\left\{
\begin{array}
[c]{l}%
dX_t  =AX_t dt+B F(t,X_t)+B dW_t
,\text{ \ \ \ }t\in\left[ 0,T\right] \\
X_0 =x,
\end{array}
\right.  \label{ornstein intro}%
\end{equation}
that is, at least formally, 
$$
(\call_t f)(x)=\frac{1}{2}(Tr Q \nabla^2 f)(x)+\<Ax+F(t,x),\nabla f(x)\>.
$$
For the sake of simplicity the paper focuses on the case of an Ornstein-Uhlenbeck process, that is $F=0$
in equation (\ref{ornstein intro}), and hints on how to handle the case of a perturbed Ornstein-Uhlenbeck process are given
throuhout the paper when necessary. From now on also in the introduction we assume $F=0$.
We also stress the fact that in this paper we aim to reduce the technical difficulties: with some efforts the case of Kolmogorov equations with
non linear term given by $\tilde \psi(t,x,\nabla v(t,x)B)$ instead of $\psi(\nabla v(t,x)B)+l(t,x)$
can be treated, but in the present paper we study only the special case of semilinear Kolmogorov equations like (\ref{Kolmo intro}).

Second order differential equations on Hilbert spaces have been extensively studied:
see e.g. the monograph \cite{DP3}. One of the main motivations for this study
in the non linear case is the connection with control theory, namely the fact
that in many cases the value
function of a stochastic optimal control problem is a solution to such an
equation.

\noindent Solutions of semilinear Kolmogororv equations
(\ref{Kolmo intro}) are studied in the literature both by an
analytic approach and by a purely probabilistic approach.
In the first direction we mention the paper \cite{Go1}, where the main
assumption is a regularizing property for the transition semigroup $P_t $, namely the strong Feller property.
In the same direction we cite also the papers
\cite{Mas} and \cite{Mas1} where equation with the special structure of (\ref{Kolmo intro}) is considered in an Hilbert space and in a Banach space respectively,
by requiring on the transition semigroup a regularizing property strictly weaker than the Strong Feller one.

\noindent For what concerns the purely probabilistic approach, semilinear Kolmogorov equations in $H$
with the special structure of (\ref{Kolmo intro}) and with non
constant $B$, are treated in the
paper \cite{fute} by means of backward stochastic
differential equations (BSDEs in the following). No regularizing assumption on the transition semigroup is imposed,
on the contrary $\psi, \,l$ and $\phi$ are assumed differentiable. The papers \cite{Mas-Ban} and \cite{Mas-SAP} are the extension of \cite{fute}
to the Banach space framework, with some restrictions on $B$, which is asked to be constant.
We notice that \cite{fute}
is the infinite dimensional extension of results in \cite{PaPe}.

In the present paper we use both the probabilistic and the analytic approach, since we want to treat
 the case of $\psi, \,l(t,\cdot)$ and $\phi$ only Lipschitz continuous
by requiring a directional regularizinging property on the transition semigroup $P_t$, similar but
weaker than the one in \cite{Mas} and \cite{Mas1}. The BSDE approach in the case of lipschitz continuous data has been used in \cite{futeGradGen}, in that paper the solution of the HJB equation is given in a sense weaker than the mild sense, since the directional derivative is addressed as a generalized gradient. So the present paper improves the results in \cite{futeGradGen} in the situations when the following assumptions in (\ref{ornstein inclusione:stimaIntro}) and
(\ref{ipotesi H-intro}) are satisfied.

Coming into more details, we assume that it holds true
\begin{equation}\label{ornstein inclusione:stimaIntro}
\left\|  Q_{t}^{-1/2}e^{tA}B\right\|  \leq c(t)\text{, for }0<t\leq T
\end{equation}
where
\[
 Q_t=\int_0^te^{sA}BB^*e^{sA^*}ds
\]
As a consequence we have a regularizing property for the transition semigroup $P_t$: it maps bounded and continuous functions
into $B$-G\^ateaux differentiable functions and for every bounded and continuous function $\phi$
\begin{equation}
\left|  \nabla^{B}P_{\tau}\left[  \phi\right]
\left(  x\right)  \xi\right|  \leq c(\tau)
\left\|  \phi\right\|  _{\infty}\left|  \xi\right|  .\label{ipotesi H-intro}%
\end{equation}The novelty of this paper towards \cite{Mas} and \cite{Mas1} is the fact that $c$,
that it is expected to blow up as $\tau$ goes to $0$, may be not locally integrable at $0$.

In order to prove existence 
and uniqueness of a mild solution $v$ of equation (\ref{Kolmo intro}),
we use the fact that if $l$, $\psi$ and $\phi$ are G\^ateaux differentiable, $v$ can be represented in terms of the solution of a suitable
forward-backward system (FBSDE in the following):
\begin{equation}\label{fbsde intro}
    \left\{\begin{array}{l}\dis dX_\tau =
AX_\tau d\tau+ B dW_\tau,\quad \tau\in
[t,T]\subset [0,T],
\\\dis
X_t=x,
\\\dis
 dY_\tau=-\psi(Z_\tau)\;d\tau-l(\tau,X_\tau)\;d\tau+Z_\tau\;dW_\tau,
  \\\dis
  Y_T=\phi(X_T).
\end{array}\right.
\end{equation}
It is well known, see e.g. \cite{PaPe} for the finite dimensional case and \cite{fute} and \cite{Mas-Ban}
for the generalization to the infinite dimensional case, respectively in the Hilbert space framework and in the Banach space framework,
that, letting $v$ be solution of the Kolmogorov equation (\ref{Kolmo intro}) when all the data $l,\,\psi$ and $\phi$ are differentiable, $v(t,x)=Y_t^{t,x}$ and $\nabla v(t,x)B=Z_t^{t,x}$.
This identification has been extended in \cite{Mas}
to the case of data continuous and with the transition semigroup satisfying the regularizing property
stated in (\ref{ipotesi H-intro}), with $c(\cdot)$
locally integrable.
These identifications are here extended, both in Hilbert and in the Banach space case, to the case of Lipschitz continuous coefficients
with the transition semigroup satisfying the regularizing property in (\ref{ipotesi H-intro}), with $c(\cdot)$ not
locally integrable.

The model we have in mind is a stocahstic heat equation on the space interval $[0,1]$, and with noise on a subdomain, whose closure is strictly contained in $[0,1]$, see Section \ref{sezionesde} for more details. This equation can be reformulated as an evolution equation both in
the Hilbert space $H$ of square integrable functions on $[0,1]$ and in the Banach space of continuous functions
on $[0,1]$.
We notice that for such a stochastic partial differential equations, suitably reformulated in an infinite dimensional 
space of functions, conditions (\ref{ornstein inclusione:stimaIntro}) and (\ref{ipotesi H-intro}) hold true, with $c(t)\sim e^{\frac{1}{t}}$,
as $t\rightarrow 0$, see the papers \cite{FZ}, \cite{Zu} and \cite{Zu1}.

The paper is organized as follows:
in Section \ref{sezionesde} some results on the Ornstein-Uhlenbeck process in $H$
are collected and the stochastic heat equation we can treat is presented, in Section \ref{applic contr 1} we present the optimal control problem we can treat in the Hilbert space framework, in Section \ref{sezioneKolmo} the Kolmogorov equation (\ref{Kolmo intro}) is solved in the Hilbert space $H$.
The we turn to the Banach space case: the choice of presenting first the Hilbert space case and then turn to the Banach space case is due to our willing of presenting the results in the simplest context, and then to pass to more complicated situations. In Section \ref{sezioneKolmoBanach} we present the Ornstein-Uhlenbeck process
in the Banach space $E$ and we solve in $E$ the Kolmogorov equation (\ref{Kolmo intro}), finally in Section \ref{applic contr 2} we 
solve the optimal control problem in the Banach space $E$: we notice that we perfom the fundamental relation in a Banach space without differentiability assumptions on the costs and on the Hamiltonian, and this is an improvement torwards
\cite{Mas-Ban}. Finally we show how our results apply to a controlled stochastic heat equation in $E$
with control and noise on a subdomain.

\section{Preliminary results on the forward equation and its semigroup}
\label{sezionesde}
We consider an Ornstein-Uhlenbeck process in a real and separable Hilbert space
$H$, that is a Markov process $X$ (also denoted $X^{t,x}$ to stress the dependence
on the initial conditions) solution to equation%
\begin{equation}
\left\{
\begin{array}
[c]{l}%
dX_\tau  =AX_\tau d\tau+BdW_\tau
,\text{ \ \ \ }\tau\in\left[  t,T\right] \\
X_t =x,
\end{array}
\right.  \label{ornstein}%
\end{equation}
where $A$ is the generator of a strongly continuous semigroup in $H$ 
and $B$ is a linear bounded operator from $\Xi$ to $H$.
\noindent We define a positive and symmetric operator%
\[
Q_{\sigma}=\int_{0}^{\sigma}e^{sA}BB^{\ast}e^{sA^{\ast}}ds.
\]
Throughout the paper we assume the following.

\begin{hypothesis}
\label{ip su AB}
\begin{enumerate}
\item  The linear operator $A$ is the generator of a strongly continuous
semigroup $\left(  e^{t A},t\geq0\right)  $ in the Hilbert space $H.$
It is well known that there exist $M>0$ and $\omega\in\mathbb{R}$ such that
$\left\Vert e^{tA}\right\Vert _{L\left(  H,H\right)  }\leq Me^{\omega t}$, for
all $t\geq0$. In the following, we always consider $M\geq1$ and $\omega\geq 0$.
\item $B$ is a bounded linear operator from $\Xi$ to $H$ and $Q_{\sigma}$ is
of trace class for every $\sigma\geq0$.
\end{enumerate}
\end{hypothesis}
The process $X^{t,x}$ is clearly time-homogeneous, and for
$0\leq t\leq \tau \leq T$ we denote by $P_{\tau-t}=P_{t,\tau}$
its transition semigroug, where for every bounded and continuous function $\phi:H\rightarrow\R$
\[
 P_{t,\tau}[\phi](x)=\E\phi(X_\tau^{t,x}).
\]
It is well known that the
Ornstein-Uhlenbeck semigroup can be represented as
\[
P_{\tau}\left[  \phi\right]  \left(  x\right)  :=\int_{H}\phi\left(  y\right)
\mathcal{N}\left(  e^{\tau A}x,Q_{\tau}\right)  \left(  dy\right), \quad \tau >0  ,
\]
and $\mathcal{N}\left(  e^{\tau A}x,Q_{\tau}\right)  \left(  dy\right)  $ denotes a
Gaussian measure with mean $e^{\tau A}x,$ and covariance operator $Q_{\tau}$.

We briefly introduce the notion of $B$-differentiability,
see e.g. \cite{Mas}. We recall that for a continuous function
$f:H\rightarrow\mathbb{R}$ the $B$-directional derivative $\nabla^{B}$ at a
point $x\in H$ in direction$\ \xi\in H$ is defined as follows:%
\[
\nabla^{B}f\left(  x;\xi\right)  =\lim_{s\rightarrow0}\frac{f\left(
x+sB\xi\right)  -f\left(  x\right)  }{s},\text{ }s\in\mathbb{R}\text{.}%
\]
A continuous function $f$ is $B$-G\^ateaux differentiable at a point $x\in H$ if
$f$ admits the $B$-directional derivative $\nabla^{B}f\left(  x;\xi\right)  $
in every directions $\xi\in \Xi$ and there exists a functional, the
$B-$gradient $\nabla^{B}f\left(  x\right)  \in\Xi^{\ast}$ such that
$\nabla^{B}f\left(  x;\xi\right)  =\nabla^Bf\left(  x\right)  \xi$.

Throughout the paper we assume the following:
\begin{hypothesis}\label{ipH su fi}
The operators $A$ and $B$ are such that
\begin{equation} 
\operatorname{Im}e^{tA}B\subset\operatorname{Im}Q_{t}^{1/2}.
\label{ornstein inclusione}
\end{equation}
As an immediate consequence it turns out the operator $Q_{t}^{-1/2}e^{tA}B$
is well defined. Assume that there exists $c:(0,T]\rightarrow\R$, such that
 $c$ is not integrable in $0$ and $c$ is bounded on any interval $I\in[0,T]$, such that
$0\notin \bar I$, where $\bar I$ is the closure of $I$, moreover we ask that $c(\cdot)$ is monotone not increasing in the interval $(0,T]$.
Assume that it holds true
\begin{equation}\label{ornstein inclusione:stima}
\left\|  Q_{t}^{-1/2}e^{tA}B\right\|  \leq c(t)\text{, for }0<t\leq T.
\end{equation}
\end{hypothesis}
As a consequence we have a regularizing property for the transition semigroup $P_t$: it maps bounded and continuous functions to $B$-G\^ateaux differentiable functions. In the following we will refer to this regularizing property as $B$-regularizing property.
\begin{lemma}Let $A$ and $B$ satisfy hypothesis \ref{ip su AB}.
\label{ipH su fi-cons} For some $c:(0,T]\rightarrow\mathbb R$
measurable and for
every $\phi\in C_{b}\left(  H\right)  $, the function $P_{\tau}\left[
\phi\right]  \left(  x\right)  $ is $B  $-differentiable
with respect to $x$, for every $0\leq t <\tau < T$ and for
every $\xi\in\Xi,$ and for $0\leq t<\tau\leq T$,%
\begin{equation}
\left|  \nabla^{B}P_{\tau}\left[  \phi\right]
\left(  x\right)  \xi\right|  \leq c(\tau)
\left\|  \phi\right\|  _{\infty}\left|  \xi\right|  .\label{ipotesi H}%
\end{equation}
then hypothesis \ref{ipH su fi-cons} is satisfied.
\end{lemma}
\dim
The proof goes on like the proof of lemma 3.4 in \cite{Mas}, where $c(t)=t^{-\alpha}$, $0<\alpha<1$.
 Note that in the present
paper we are mainly concerned with the case $c(\cdot)\notin L^1([0,T])$, but this aspect does not enter the proof of the present result. The case
$c(\cdot)\in L^1([0,T])$ is treatable as in \cite{Mas}.
Note also that the $B$-regularizing property for the transition semigroup $P_t$ stated in (\ref{ipotesi H})
is equivalent to the inclusion in (\ref{ornstein inclusione}).
\qed

\begin{remark}\label{remark:OUpertReg}
In the case of $A$ and $B$ satisfying hypothesis \ref{ipH su fi}
the transition semigroup of the perturbed Ornstein Uhlenbeck process
\begin{equation}\label{ornstein-pert}
\left\{
\begin{array}
[c]{l}%
dX_\tau  =AX_\tau d\tau+BF(\tau,X_\tau)+BdW_\tau
,\text{ \ \ \ }\tau\in\left[  t,T\right] \\
X_t =x,
\end{array}
\right.
\end{equation}
satisfies the $B$-regularizing property proved in lemma \ref{ipH su fi-cons},
and stated in (\ref{ipotesi H}), by assuming that
$F$ is jointly continuous in $t$ and $x$ and lipschitz continuous in $x$ uniformly
with respect to $t$ and G\^ateaux differentiable with respect to $x$. For
the proof see \cite{Mas1}, Section 4, Theorem 4.3.
\end{remark}

The model we have in mind is a semilinear heat equation.
Namely let $\calo$ be a subinterval of the interval $[0,1]$ such that 
$\bar \calo_0 \subsetneq [0,1]$. In the following $\calo=[a,b],\, 0<a<b<1$.
We denote by $H$ the Hilbert space $L^2([0,1])$ and the equation
\begin{equation}\label{heat equation}
 \left\{
  \begin{array}{l}
  \dis
\frac{ \partial y}{\partial s}(s,\xi)= \Delta y(s,\xi)+ 1_{\calo}(\xi)f(s,y(s,\xi))
+ 1_{\calo}(\xi)\frac{ \partial W
}{\partial s}(s,\xi), \qquad s\in [t,T],\;
\xi\in [0,1],
\\\dis
y(t,\xi)=x(\xi),
\\\dis
\dfrac{\partial}{\partial\xi} y(s,\xi)=0, \quad \xi=0,\;\xi=1.
\end{array}
\right.
\end{equation}
Here $\frac{ \partial W
}{\partial s}(s,\xi)$ is a space time white noise and $f:[0,T]\times [0,1]\rightarrow \R$ is a continuous
function such that $f(s,\cdot):[0,1]\rightarrow \R$ is differentiable with derivative uniformly bounded with respect to $s\in[0,T]$.
Let $F$ be the evaluation operator associated to $f$ and $B$ the multiplication
operator associated to $1_{\calo}$: for every $h\in H$, $Bh(\xi)=1_{\calo}(\xi)h(\xi)$.
With this definition of $F$ and $B$, equation (\ref{heat equation}) can be written in an abstract way in $H$ as
a perturbed Ornstein Uhlenbeck process, see equation (\ref{ornstein-pert}),
where $A$ is the Laplace operator with Neumann boundary conditions and
$W$ is a cylindrical Wiener process in $H$.

\section{The optimal control problem}
\label{applic contr 1}

We consider the following controlled state equation
\begin{equation}
\left\{
\begin{array}
[c]{l}%
dX^{u}_\tau  =\left[  AX^{u}_\tau +B u_\tau
 \right]  d\tau+BdW_\tau ,\text{ \ \ \ }\tau\in\left[  t,T\right] \\
X^{u}_t  =x.
\end{array}
\right.  \label{sdecontrolforte}%
\end{equation}
The solution of this equation will be denoted by
$X_\tau^{u,t,x}$ or simply by $X^{u}_\tau$. $X$ is also called
the state, $T>0,$ $t\in\left[  0,T\right]$ are fixed.
The process $u$ represents the control and it is an $\left(\mathcal{F}_{\tau}\right)_{\tau}$-predictable
process with values in a closed and bounded set $U$ of the Hilbert space $\Xi$, such that $\vert u\vert \leq R$.
The occurrence of the operator
$B$ in the control term is imposed by our techniques
and, among other facts, it allows to study
the optimal control problem related by means of BSDEs.

Beside equation (\ref{sdecontrolforte}),\ define the cost
\begin{equation}
J\left(  t,x,u\right)  =\mathbb{E}\int_{t}^{T}[
l\left(s,X^{u}_s\right)+g(u_s) ]ds+\mathbb{E}\phi\left(X^{u}_T\right). 
\label{cost}%
\end{equation}
for real functions $l$ on $[0,T]\times H$, $\mathbb{\phi}$ on $H$ and $g$ on $U$.
The control problem in strong formulation is to
minimize this functional $J$ over all admissible controls $u$.
We make the following assumptions on the cost $J$.

\begin{hypothesis}
\label{ip costo}

\begin{enumerate}
\item  The function $\mathbb{\phi}:H\rightarrow\mathbb{R}$ is
bounded and lipschitz continuous;

\item $l:[0,T]\times H\rightarrow\mathbb{R}$ is bounded,
lipschitz continuous with respect to $x\in H$
uniformly with respect to $t\in[0,T]$;

\item $g: U\rightarrow\mathbb{R}$ is bounded and continuous.
\end{enumerate}
\end{hypothesis}
We deduce immediately  that equation (\ref{sdecontrolforte}) admits a unique mild solution, for
every admissible control $u$.

\noindent We denote by $J^{\ast}\left(  t,x\right)  =\inf_{u\in\mathcal{A}%
_{d}}J\left(  t,x,u\right)  $ the value function of the problem and, if it
exists, by $u^{\ast}$ the control realizing the infimum, which is called
optimal control.

We define in a classical way the Hamiltonian function relative to the above
problem:%
\begin{equation}\label{hamilton}
\psi\left(z\right)  =\inf_{u\in U}\left\{  g\left(u\right)
+z u\right\}\quad \forall z\in \Xi.
\end{equation}
By our assumptions the Hamiltonian function is lipschitz continuous. 

We define
\begin{equation}\label{defdigammagrande}
\Gamma(z)=\left\{ u\in U: zu+g(u)= \psi(z)\right\};
\end{equation}
if $\Gamma(z) \neq \emptyset$ for every $z\in \Xi$, by \cite{AuFr}, see Theorems 8.2.10 and
8.2.11, $\Gamma$ admits a measurable selection, i.e. there exists
a measurable function $\gamma: H \rightarrow U$ with
$\gamma(z)\in \Gamma(z)$ for every $z\in \Xi$.

\section{The semilinear Kolmogorov equation}
\label{sezioneKolmo}
The aim of this section is to present existence and uniqueness results for the solution of
the Hamilton Jacobi Bellman equation ( HJB in the following ) related to the optimal control
problem presented in section \ref{applic contr 1}.

More precisely, let $\mathcal L $ be the generator of the transition
semigroup $P_t$, that is, at least formally, 
$$
(\call f)(x)=\frac{1}{2}(Tr BB^* \nabla^2 f)(x)+\<Ax,\nabla f(x)\>.
$$
Let us consider the following equation
\begin{equation}
\left\{
\begin{array}
[c]{l}%
-\frac{\partial v}{\partial t}(t,x)=\mathcal{L} v\left(  t,x\right)
+\psi\left(  \nabla^{B}v\left(  t,x\right)  \right)+l(t,x)  ,\text{ \ \ \ \ }t\in\left[  0,T\right]
,\text{ }x\in H\\
v(T,x)=\phi\left(  x\right)  ,
\end{array}
\right.  \label{Kolmo}%
\end{equation}
We introduce the notion of mild solution of the non linear Kolmogorov
equation (\ref{Kolmo}), see e.g. \cite{fute} and also \cite{Mas} for the definition of
mild solution when $\psi$ depends only on $\nabla^{B}v$ and not on $\nabla v$. 
Since $\mathcal{L}$ is (formally) the generator of
$P_{t}$, the variation of constants formula for (\ref{Kolmo}) is:%
\begin{equation}
v(t,x)=P_{t,T}\left[  \phi\right]  \left(  x\right)  +\int_{t}^{T}%
P_{t,s}\left[  \psi(\nabla^{B }v\left(  s,\cdot\right))  \right]  \left(  x\right)  ds
+\int_{t}^{T}
P_{t,s}\left[  l\right(s,\cdot)]  \left(  x\right)  ds,\text{\ \ }t\in\left[  0,T\right]  ,\text{ }x\in H. \label{solmildkolmo}%
\end{equation}
We use this formula to give the notion of mild
solution for the non linear Kolmogorov equation (\ref{Kolmo}); we have also to
introduce some spaces of continuous functions, where we seek the
solution of (\ref{Kolmo}).

\noindent As stated before we focus on the case $c(\cdot )\notin L^1([0,T])$:
let $C_{c(\cdot)}\left(  \left[  0,T\right]  \times
H\right)  $ be the linear space of continuous functions $f:\left[  0,T\right)
\times H\rightarrow\mathbb{R}$ such that 
$$
\sup_{t\in\left[  0,T\right]  }%
\sup_{x\in H}(c\left(  T-t\right))  ^{-1}\left\vert f\left(  t,x\right)
\right\vert <+\infty.$$
$C_{c(\cdot)}\left(  \left[  0,T\right]  \times H\right)$
endowed with the norm
\[
\left\Vert f\right\Vert _{C_{c(\cdot)}}=\sup_{t\in\left[  0,T\right]  }%
\sup_{x\in H}(c\left(  T-t\right))  ^{-1}\left\vert f\left(  t,x\right)
\right\vert ,
\]
is a Banach space.

\noindent We consider also the linear space $C_{c(\cdot)}^{s}\left(  \left[
0,T\right]  \times H,\Xi^{\ast}\right)  $ of the mappings $L:\left[
0,T\right)  \times H\rightarrow \Xi^{\ast}$ such that for every $\xi\in \Xi$,
$L\left(  \cdot,\cdot\right)  \xi\in C_{c(\cdot)}\left(  \left[  0,T\right]
\times H\right)  $. The space $C_{c(\cdot)}^{s}\left(  \left[  0,T\right]
\times H,\Xi^{\ast}\right)  $ turns out to be a Banach space if it is endowed
with the norm
\[
\left\Vert L\right\Vert _{C_{c(\cdot)}\left(  \Xi^{\ast}\right)  }=\sup
_{t\in\left[  0,T\right]  }\sup_{x\in H}(c\left(  T-t\right))  ^{-1}
\left\Vert L\left(  t,x\right)  \right\Vert _{\Xi^{\ast}}.
\]
In other words, $C_{c(\cdot)}^{s}\left(  \left[  0,T\right]  \times H,\Xi^{\ast
}\right)  $ can be identified with the space of the operators

\noindent$L\left( H,C_{c(\cdot)}\left(  \left[  0,T\right]  \times H\right)
\right)  $.

\begin{definition}
\label{defsolmildkolmo}
%Let $c(\cdot)\in L^1\left[  0,T\right)  $.
We say that a
function $v:\left[  0,T\right]  \times H\rightarrow\mathbb{R}$ is a mild
solution of the non linear Kolmogorov equation (\ref{Kolmo}) if the following
are satisfied:

\begin{enumerate}
\item $v\in C_{b}\left(  \left[  0,T\right]  \times H\right)  $;

\item $\nabla^{B}v\in C_{c(\cdot)}^{s}\left(  \left[
0,T\right]  \times H,\Xi^{\ast}\right)  $: in particular this means that for
every $t\in\left[  0,T\right)  $, $v\left(  t,\cdot\right)$ is $B$-differentiable;

\item  equality (\ref{solmildkolmo}) holds.
\end{enumerate}
\end{definition}

Existence and uniqueness of a mild solution of equation (\ref{Kolmo})
are related to the study of the following 
forward-backward system: for given $t\in [0,T]$ and $x\in H$,
\begin{equation}\label{fbsde}
    \left\{\begin{array}{l}\dis dX_\tau =
AX_\tau d\tau+B dW_\tau,\quad \tau\in
[t,T]\subset [0,T],
\\\dis
X_t=x,
\\\dis
 dY_\tau=-\psi(Z_\tau)\;d\tau-l(\tau,X_\tau)\;d\tau+Z_\tau\;dW_\tau,
  \\\dis
  Y_T=\phi(X_T),
\end{array}\right.
\end{equation}
and to the identification of $Z_t^{t,x}=\nabla_x Y_t^{t,x} B $.
We extend the definition of $X$ setting
$X_s=x$ for $0\leq s\leq t$. The second equation in
(\ref{fbsde}), namely
\begin{equation}\label{bsde}
    \left\{\begin{array}{l}\dis
 dY_\tau=-\psi(Z_\tau)\;d\tau-l(\tau,X_\tau)\;d\tau+Z_\tau\;dW_\tau,
 \qquad \tau\in [0,T],
  \\\dis
  Y_T=\phi(X_T),
\end{array}\right.
\end{equation}
is of backward type.
Under suitable assumptions on the coefficients
 $\psi:\Xi
\rightarrow\mathbb{R}$, $l:[0,T]\times H\rightarrow\mathbb{R}$
and $\mathbb{\phi}:H\rightarrow\mathbb{R}$
we will look for a solution consisting of a pair of predictable processes,
taking values in $\mathbb{R}\times H$, such that $Y$ has
continuous paths and
\[
\|\left( Y,Z\right)\|^2_{\mathbb{K}_{cont}}:=
\mathbb{E}\sup_{\tau\in\left[ 0,T\right] }\left\vert Y_{\tau}\right\vert
^{2}+\mathbb{E}\int_{0}^{T}\left\vert Z_{\tau}\right\vert ^{2}d\tau<\infty,
\]
see e.g. \cite{PaPe1}.
In the following we denote by $\mathbb{K}_{cont}\left(
\left[ 0,T\right] \right)$ the space of such processes.

The solution of (\ref{fbsde}) will be denoted by $(X_\tau, Y_\tau, Z_\tau)_{\tau\in[0,T]}$, or,
to stress the dependence on the initial time $t$ and on the
initial datum $x$, by $(X_\tau^{t,x}, Y_\tau^{t,x}, Z_\tau^{t,x})_{\tau\in[0,T]}$.
In the following we refer to \cite{fute} for the definition of the class 
$\calg(H)$ of G\^ateaux differentiable functions $f:H\rightarrow \R$
with strongly continuous derivative. We make differentiability assumptions on the coefficients that we are going to remove in the sequel.
\begin{hypothesis}
\label{ip diffle agg}The map $\psi:\Xi
\rightarrow\mathbb{R}$ is in $ \calg(\Xi)$, and the maps $l(t,\cdot): H\rightarrow\mathbb{R}$ and $\phi:H\rightarrow\mathbb{R}$
belong to $ \calg(H)$.
\end{hypothesis}
Note that since the Hamiltonian function is lipschitz continuous
and by lipschitz assumptions on $l$ and $\phi$, see hypothesis \ref{ip costo}, it turns out that 
$\nabla \psi(z)$, $\nabla_xl(t,x)$ and $\nabla\phi(x)$  are bounded.

\noindent It is well known, see e.g. \cite{PaPe} and \cite{fute} for the infinite dimensional extension,
that under hypothesis \ref{ip costo}
the BSDE (\ref{bsde}) admits a unique solution $(Y_\tau^{t,x},Z_\tau^{t,x})\in \mathbb{K}_{cont}\left(
\left[ 0,T\right] \right)$ and, if we further assume \ref{ip diffle agg},
setting $v(t,x):=Y_t^{t,x}$,
it turns out that $v$ is the unique mild solution of equation (\ref{Kolmo}), 
and $\nabla^{B}v(t,x)=Z_t^{t,x}$. This leads immediately to the following lemma:
\begin{lemma}\label{lemmaKolmo-infsup}
 Assume that hypotheses \ref{ip su AB} and \ref{ip costo} hold true.
Let $\psi_n$, $\phi_n$ and $l_n(\tau, \cdot)$ the inf-sup convolution of
$\psi$, $\phi$ and $l$ respectively. Then the Kolmogorov equation
\begin{equation}
\left\{
\begin{array}
[c]{l}%
-\frac{\partial v}{\partial t}(t,x)=\mathcal{L} v\left(  t,x\right)
+\psi_n\left(  \nabla^{B}v\left(  t,x\right)  \right)+l_n(t,x)  ,\text{ \ \ \ \ }t\in\left[  0,T\right]
,\text{ }x\in H\\
v(T,x)=\phi_n\left(  x\right)  ,
\end{array}
\right.  \label{Kolmo_n}%
\end{equation}
admits a unique mild solution, with bounded G\^ateaux derivative.
\end{lemma}
\dim
 We recall, see e.g. \cite{LL} 
and \cite{DP3}, that the inf-sup convolution
$\phi_{n}$, $\psi_n$ and $l_n$ are defined respectively by 
\begin{equation}
\phi_{n}\left(  x\right)  =\sup_{z\in
H}\left\{  \inf_{y\in H}\left[  \phi\left(  y\right)  +n\frac{\left|
z-y\right|  _{H}^{2}}{2}\right]  -n\left|  x-z\right|  _{H}^{2}
\right\}  , \label{infsupconv-phi}%
\end{equation}
by
\begin{equation}
\psi_{n}\left(  z\right)  =\sup_{x\in\Xi}\left\{  \inf_{y\in \Xi}\left[  \psi\left(  y\right)  +n\frac{\left|
x-y\right|  _{\Xi}^{2}}{2}\right]  -n\left|  x-z\right|  _{\Xi}^{2}
\right\}  , \label{infsupconv-psi}%
\end{equation}
and by
\begin{equation}
l_{n}\left(  t,x\right)  =\sup_{z\in
H}\left\{  \inf_{y\in H}\left[  l\left(t,  y\right)  +n\frac{\left|
z-y\right|  _{H}^{2}}{2}\right]  -n\left|  x-z\right|  _{H}^{2}
\right\}  , \label{infsupconv-l}%
\end{equation}
for every $t\in[0,T]$.
It turns out that $\phi_{n},\,l_n(t,\cdot)\in\calg(H)$ and $\psi_n\in\calg(\Xi)$,
and so by \cite{fute},
there exists a unique solution to the forward backward system
\begin{equation}\label{fbsde_n}
    \left\{\begin{array}{l}\dis dX_\tau =
AX_\tau d\tau+BF(\tau,X_\tau) d\tau+ B dW_\tau,\quad \tau\in
[t,T]\subset [0,T],
\\\dis
X_t=x,
\\\dis
 dY^n_\tau=-\psi_n(Z^n_\tau)\;d\tau-l_n(\tau,X_\tau)\;d\tau+Z^n_\tau\;dW_\tau,
  \\\dis
  Y^n_T=\phi_n(X_T),
\end{array}\right.
\end{equation}
that we denote by $(X^{t,x},Y^{n,t,x},Z^{n,t,x})$ and by setting
$v_n(t,x)=Y^{n,t,x}_t$, it turns out that $v_n$ is the unique mild solution to
the Kolmogorov equation (\ref{Kolmo_n}).
\qed

Next we prove that equation (\ref{Kolmo}) admits a unique mild solution, in the case when the final conditions is lipschitz continuous, $\psi$ is assumed further to be differentiable and $l=0$, and that this solution is given in terms of the solution of the related FBSDE (\ref{fbsde}).

\begin{theorem}
\label{teoKolmo}
 Assume that hypotheses \ref{ip su AB} and \ref{ip costo} hold true, assume that in equation (\ref{Kolmo})
$l=0$ and $\psi$ is G\^ateaux differentiable.
Then equation (\ref{Kolmo}) admits a unique mild solution $v(t,x)$ according to definition
(\ref{defsolmildkolmo}), and $v(t,x)$ is given by
\[
 v(t,x):=Y^{t,x}_t
\]
with $Y^{t,x}$ solution to the BSDE in (\ref{fbsde}).
\end{theorem}
\dim
We build the solution $v$ by an approximation procedure.
We let $(\phi_n)_n$ be the inf-sup convolutions of $\phi$, see (\ref{infsupconv-phi}),
and we consider the Kolmogorov equation
\begin{equation}
\left\{
\begin{array}
[c]{l}%
-\frac{\partial v}{\partial t}(t,x)=\mathcal{L} v\left(  t,x\right)
+\psi\left(  \nabla^{B}v\left(  t,x\right)  \right) ,\text{ \ \ \ \ }t\in\left[  0,T\right]
,\text{ }x\in H\\
v(T,x)=\phi_n\left(  x\right)  ,
\end{array}
\right.  \label{Kolmo_phi_n}%
\end{equation}
whose mild solution $v_n$ by lemma \ref{lemmaKolmo-infsup} exists and is unique, and it is given by
\[
 v_n(t,x):=Y^{n,t,x}_t
\]
where $(X^{t,x},Y^{n,t,x},Z^{n,t,x})$ solve the FBSDE
\begin{equation}\label{fbsde_phi_n}
    \left\{\begin{array}{l}\dis dX_\tau =
AX_\tau d\tau+ B dW_\tau,\quad \tau\in
[t,T]\subset [0,T],
\\\dis
X_t=x,
\\\dis
 dY^{n,t,x}_\tau=-\psi(Z^{n,t,x}_\tau)\;d\tau+Z^{n,t,x}_\tau\;dW_\tau,
  \\\dis
  Y^{n,t,x}_T=\phi_n(X_T).
\end{array}\right.
\end{equation}
Let us set, for $h\in\Xi$,  $F^{n,t,x,h}_\cdot:=<\nabla_x Y^{n,t,x}_\cdot ,Bh>$ and 
$V^{n,t,x,h}_\cdot:=<\nabla_x Z^{n,t,x}_\cdot,Bh>$.
The pair of processes $(F^{n,t,x,h},V^{n,t,x,h})$ solve the following BSDE
\begin{equation}\label{fbsde_phi_n-diffle.}
    \left\{\begin{array}{l}
 dF^{n,t,x,h}_\tau=-\nabla\psi(Z^{n,t,x,h}_\tau)V^{n,t,x,h}_\tau\;d\tau+V^{n,t,x,h}_\tau\;dW_\tau,
  \\\dis
  F^{n,t,x,h}_T=\nabla\phi_n(X_T^{t,x})e^{(T-t)A}Bh,
\end{array}\right.
\end{equation}
It is standard to see that
\[
 Y^{n,t,x}\rightarrow Y^{t,x} \text{ in }L^2(\Omega,C([0,T])),\; Z^{n,t,x}\rightarrow Z^{t,x}
\text{ in }L^2(\Omega\times[0,T]).
\]
and so, by taking a subsequence $n_k$, also 
\[
 Z^{n_k,t,x}\rightarrow Z^{t,x}\, dt\times d\P
\text{ a.e. }.
\]
For what concerns $F^{n,t,x}$ and $V^{n,t,x}$, since by the lipschitzianity of $\phi$ we have
that $\nabla\phi_n$ is bounded uniformly with respect to $n$, by assumptions \ref{ip su AB}
and \ref{ip costo} on $\psi$, by standard results on BSDEs
we have that there exists a constant $C$ independent on $n$ such that 
\begin{equation}\label{stimaFnVn}
 \Vert F^{n,t,x}\Vert_{L^2(\Omega,C([0,T]))}
+ \Vert V^{n,t,x}\Vert_{L^2(\Omega\times [0,T])}\leq C.
\end{equation}
% So by the Banach-Alaoglu theorem, unless taking a subsequence $n_k$,
% there exists $V\in L^2(\Omega,\times [0,T])$ such that
% \[
%  V^{n_k,t,x}\rightarrow V^{t,x} \text{ weakly in }L^2(\Omega,\times [0,T]).
% \]
Next we investigate the convergence of $F^{n,t,x,h}$.
By writing the BSDE satisfied by $F^{n,t,x,h}-F^{j,t,x,h}\,n,j\geq 1$, integrating between $t$ and $T$ and taking expectation, we get
\begin{align*}
  F^{n,t,x,h}_t-F^{j,t,x,h}_t& = \E\left[\nabla\phi_n(X_T^{t,x})-
\nabla\phi_j(X_T^{t,x})\right]e^{(T-t)A}Bh \\
&+ \E\int_t^T \left( 
\nabla\psi(Z^{n,t,x}_\tau)V^{n,t,x,h}_\tau-
\nabla\psi(Z^{j,t,x}_\tau)V^{j,t,x,h}_\tau\right)\,d\tau
\end{align*}
and consequently
\begin{align*}
 \vert F^{n,t,x,h}_t-F^{j,t,x,h}_t\vert &\leq \vert\E\left[\nabla\phi_n(X_T^{t,x})-
\nabla\phi_j(X_T^{t,x})\right]e^{(T-t)A}Bh\vert\\
& +\vert\E \int_t^T \left( 
\nabla\psi(Z^{n,t,x}_\tau)V^{n,t,x,h}_\tau- \nabla\psi(Z^{j,t,x}_\tau)V^{j,t,x,h}_\tau\right)\,d\tau
\vert =I+II
\end{align*}
We start by estimating $I$:
\begin{align*}
 I&=\vert \E<\nabla^B\left(\phi_n(X_T^{t,x})
-\phi_j(X_T^{t,x})\right),h>\vert=\vert\nabla^B<\E\left(\phi_n-\phi_j\right)(X_T^{t,x}),h>\vert\\
&=\vert <\nabla^B P_{t,T}\left[\phi_n-\phi_j\right](x),h>\vert\leq c(T-t)\Vert\phi_n-\phi_j\Vert_\infty
\vert h\vert
\end{align*}
In order to estimate $II$, we notice that 
since for every $n\geq 1$ the final datum $\phi$ and the hamiltonian function $\psi$ are differentiable, 
then 
we have the following identication for $F^{n,t,x}_t$:
\begin{equation}\label{identifications}
 <\nabla^B Y^{n,t,x}_t,h>:=F^{n,t,x,h}_t=<\nabla^B v_n(t,x),h>= <Z^{n,t,x}_t,h>,
\end{equation}
where in mild form $v_n$ satisfies the integral equation
\[
 v_n(t,x)=P_{t,T}[\phi_n](x)+\int_t^TP_{t,s}[\psi\left(\nabla^Bv_n(s,\cdot)\right)](x)\,ds.
\]
By taking the directional derivative $\nabla^B$ we get
\[
 \nabla^Bv_n(t,x)=\nabla^B P_{t,T}[\phi_n](x)+\int_t^T \nabla^B P_{t,s}[\psi\left(\nabla^Bv_n(s,\cdot)\right)](x)\,ds
\]
So $II$ can be rewritten in different ways as
\begin{align*}
 II=&\vert\E \int_t^T \left( 
\nabla\psi(Z^{n,t,x}_\tau)V^{n,t,x,h}_\tau- \nabla\psi(Z^{j,t,x}_\tau)V^{j,t,x,h}_\tau\right)\,d\tau\vert\\
=&\vert\E \int_t^T <\nabla^B\left( 
\psi(Z^{n,t,x}_\tau)-\psi(Z^{j,t,x}_\tau)\right),h>\,d\tau\vert\\
=&\vert <\nabla^B\int_t^T \E\left( 
\psi(Z^{n,t,x}_\tau)-\psi(Z^{j,t,x}_\tau)\right)\,d\tau,h>\vert\\
=&\vert \int_t^T <\nabla^B P_{t,\tau}\left[ 
\psi(\nabla^Bv_n(t,\cdot))-\psi(\nabla^Bv_j(t,\cdot))\right],h>\,d\tau\vert
\end{align*}
We fix $0<\delta<T-t$ which will be chosen in the following and in order to estimate II we start by estimating
\[
 \vert\E \int_{t+\delta}^T \left( 
\nabla\psi(Z^{n,t,x}_\tau)V^{n,t,x,h}_\tau- \nabla\psi(Z^{j,t,x}_\tau)V^{j,t,x,h}_\tau\right)\,d\tau
\vert 
\]
 and we follow the calculations in \cite{Mas}, lemma 3.4. So
\begin{align*}
&\vert\E \int_{t+\delta}^T \left( 
\nabla\psi(Z^{n,t,x}_\tau)V^{n,t,x,h}_\tau- \nabla\psi(Z^{j,t,x}_\tau)V^{j,t,x,h}_\tau\right)\,d\tau
\vert \\
&=\vert \int_{t+\delta}^T <\nabla^B P_{t,\tau}\left[ 
\psi(\nabla^Bv_n(t,\cdot))-\psi(\nabla^Bv_j(t,\cdot))\right],h>\,d\tau\vert\\
&=\vert \int_{t+\delta}^T\int_{H}\left(\psi\left(\nabla^B v_n(s,  y+e^{(s-t)A}x)\right)  -
\psi\left(\nabla^B v_j(s,  y+e^{(s-t)A}x)\right) \right) \\
& \left\langle Q_{s-t}^{-1/2}%
e^{(s-t)A}Bh,Q_{s-t}^{-1/2}y\right\rangle \mathcal{N}\left(  0,Q_{s-t}\right)  \left(
dy\right)\,ds\vert\\
  &\leq \int_{t+\delta}^T\left(\int_{H}\vert\psi\left(\nabla^B v_n(s,  y+e^{(s-t)A}x)\right)  -
\psi\left(v_j(s,  y+e^{(s-t)A}x)\right) \vert^2 \mathcal{N}\left(  0,Q_{s-t}\right) \right)^{\frac{1}{2}}\\
 &\left(\int_{H}\vert\left\langle Q_{s-t}^{-1/2}%
e^{(s-t)A}Bh,Q_{s-t}^{-1/2}y\right\rangle\vert^2 \mathcal{N}\left(  0,Q_{s-t}\right)  \left(
dy\right)\right)^{\frac{1}{2}}\,ds
\end{align*}
where the last passage follows by the Cauchy-Schwartz inequality.
So, by the Lipschitz property of $\psi$, and by the
 estimate (\ref{ornstein inclusione:stima}), we get
\begin{align*}
&\vert\E \int_{t+\delta}^T \left( 
\nabla\psi(Z^{n,t,x}_\tau)V^{n,t,x,h}_\tau- \nabla\psi(Z^{j,t,x}_\tau)V^{j,t,x,h}_\tau\right)\,d\tau
\vert \\
&\leq C  \int_{t+\delta}^T c(s-t)\left(\int_{H}\vert\nabla^B v_n\left(s,  y+e^{(s-t)A}x)\right)  -
\nabla^B v_j(s,  y+e^{(s-t)A}x) \vert^2 \mathcal{N}\left(  0,Q_{s-t}\right) \right)^{\frac{1}{2}}\\
&\leq C c(\delta) \E\int_{t+\delta}^T \vert Z^{n,t,x}_s-Z^{j,t,x}_s\vert^2\,ds .
\end{align*}
where in the last passage we have used the identification (\ref{identifications}) of $Z^{n,t,x}$
with $\nabla^B v_n$, and the fact that $c(\cdot)$ is monotone not increasing.
Now since $Z^{n,t,x}$ is a Cauchy sequence in $L^2(\Omega\times [0,T])$, for every $\varepsilon>0$
it is possible to choose $n,j \geq \bar n$ such that
\[
 \E\int_{t}^T \vert Z^{n,t,x}_s-Z^{j,t,x}_s\vert^2\,ds\leq \frac{\varepsilon}{C c(\delta)}.
\]
We can conclude that, choosen $\varepsilon>0$
\begin{equation}\label{primastimaII}
 \vert\E \int_{t+\delta}^T \left( 
\nabla\psi(Z^{n,t,x}_\tau)V^{n,t,x,h}_\tau- \nabla\psi(Z^{j,t,x}_\tau)V^{j,t,x,h}_\tau\right)\,d\tau
\vert \leq \varepsilon,
\end{equation}
independently on the choice of $\delta$. We now estimate
\[
\vert\E \int_t^{t+\delta} \left( 
\nabla\psi(Z^{n,t,x}_\tau)V^{n,t,x,h}_\tau- \nabla\psi(Z^{j,t,x}_\tau)V^{j,t,x,h}_\tau\right)\,d\tau
\vert.
\]
We notice that by the Lipschitzianity of $\psi$ the derivative $\nabla\psi$ is
bounded, and that $V^{j,t,x,h}$ and $V^{n,t,x,h}$ are uniformly
bounded in $L^2(\Omega \times [0,T])$, so for any
$ \varepsilon >0$ we can choose $\delta$ such that
\[
 \vert\E \int_t^{t+\delta} \left( 
\nabla\psi(Z^{n,t,x}_\tau)V^{n,t,x,h}_\tau- \nabla\psi(Z^{j,t,x}_\tau)V^{j,t,x,h}_\tau\right)\,d\tau
\vert\leq \varepsilon.
\]
and so we can conclude that for $n,j\geq \bar n$
\[
 II\leq2\varepsilon.
\]
We have proved that
\begin{align*}
 \vert F^{n,t,x,h}_t-F^{j,t,x,h}_t\vert &\leq c(T-t)\Vert\phi_n-\phi_j)\Vert\vert h\vert+\varepsilon
\end{align*}
So by the identification in (\ref{identifications})
we get that $\nabla^B v_n$ is a Cauchy sequence in the space of the operators $C^s_{c(\cdot)
\left( [0,T]\times H,\Xi^*\right)}$ and
so it converges, to a limit that we denote by $L(t,x)$.
We prove that $L(t,x)$ coincides with $\nabla^B v(t,x)$.
We have, for all $t\in[0,T],\,x\in H$ and for all $h\in\Xi$
\begin{align*}
 \dfrac{v_n(t,x+sBh)-v_n(t,x)}{s}&=\int_0^1 <\nabla v_n(t, x+\alpha s Bh),Bh>\,d\alpha\\
&=\int_0^1< \nabla^B v_n(t, x+\alpha s Bh),h>\,d\alpha.
\end{align*}
Passing to the limit as $s\rightarrow 0$ on both sides we get
\[
 <\nabla^Bv_n(t,x),h>=L(t,x)h,
\]
and the proof is concluded.
\qed

\begin{remark}\label{remark:nablaBbounded}
It is immediate to see that $\nabla^B v(t,x) $ is bounded. In fact we have seen in the proof of Theorem \ref{teoKolmo} that $\nabla^Bv(t,x)$ is the pointwise limit of $\nabla^B v_n(t,x)$, that $\nabla^B v_n(t,x)=F^{n,t,x}_t$
and $F^{n,t,x}_t$ is bounded, uniformly with respect to $n$. So we deduce that that $\nabla^Bv(t,x)$
is bounded. The same holds true for the solution whose existence and uniqueness is proved in Theorem
\ref{teoKolmo:gen} under more general assumptions than in Theorem \ref{teoKolmo}.
\end{remark}

In Theorem \ref{teoKolmo} we have proved that the Kolmogorov equation admits a unique mild solution,
assuming the running cost $l=0$ and $\psi$ G\^ateaux differentiable.
The result is still true in the general context of the present paper, that is removing these additional assumptions.
\begin{theorem}
\label{teoKolmo:gen}
 Assume that hypotheses \ref{ip su AB} and \ref{ip costo} hold true.
Then equation (\ref{Kolmo}) admits a unique mild solution $v(t,x)$ according to definition
(\ref{defsolmildkolmo}), and $v(t,x)$ is given by
\[
 v(t,x):=Y^{t,x}_t
\]
with $Y^{t,x}$ solution to the BSDE in (\ref{fbsde}).
\end{theorem}
\dim The case of $l\neq 0$ and lipschitz continuous with respect to $x$, uniformly with respect to $t$, can be handled as we have treated the final datum $\phi$. So we take $l=0$ and we look at removing the differentiability assumption on the Hamiltonian function $\psi$. We approximate $\phi$ and $\psi$ with their
inf-sup convolutions, defined respectively in (\ref{infsupconv-phi}) and in (\ref{infsupconv-psi}).
By Lemma \ref{lemmaKolmo-infsup}, the approximating Kolmogorov equation
\[
 \left\{
\begin{array}
[c]{l}%
-\frac{\partial v}{\partial t}(t,x)=\mathcal{L} v\left(  t,x\right)
+\psi_n\left(  \nabla^{B}v\left(  t,x\right)  \right)+l(t,x)  ,\text{ \ \ \ \ }t\in\left[  0,T\right]
,\text{ }x\in H\\
v(T,x)=\phi_n\left(  x\right)  ,
\end{array}
\right.
\]
admits a unique mild solution $v_n$ satisfying
 \begin{equation*}
v_n(t,x)=(P_{t,T}\left[  \phi_n\right]  )\left(  x\right)  
+\int_t^T P_{t,s}[\psi_n\left(\nabla^Bv_n(s,\cdot)\right)](x)\,ds.
 \end{equation*}
Again by Lemma \ref{lemmaKolmo-infsup} $v_n$ is given by
\[
 v_n(t,x):=Y^{n,t,x}_t
\]
where $(X^{t,x},Y^{n,t,x},Z^{n,t,x})$ solve the FBSDE
\begin{equation}\label{fbsde_phi_n-psi_n}
    \left\{\begin{array}{l}\dis dX_\tau =
AX_\tau d\tau+B dW_\tau,\quad \tau\in
[t,T]\subset [0,T],
\\\dis
X_t=x,
\\\dis
 dY^{n,t,x}_\tau=-\psi_n(Z^{n,t,x}_\tau)\;d\tau+Z^{n,t,x}_\tau\;dW_\tau,
  \\\dis
  Y^{n,t,x}_T=\phi_n(X_T).
\end{array}\right.
\end{equation}
Let us set, for $h\in\Xi$,  $F^{n,t,x,h}_\cdot:=<\nabla_x Y^{n,t,x}_\cdot ,Bh>$ and 
$V^{n,t,x,h}_\cdot:=<\nabla_x Z^{n,t,x}_\cdot,Bh$.
The pair of processes $(F^{n,t,x,h},V^{n,t,x,h})$ solve the following BSDE
\begin{equation}\label{fbsde_phi_npsi_n-diffle.}
    \left\{\begin{array}{l}
 dF^{n,t,x,h}_\tau=-\nabla\psi_n(Z^{n,t,x,h}_\tau)V^{n,t,x,h}_\tau\;d\tau+V^{n,t,x,h}_\tau\;dW_\tau,
  \\\dis
  F^{n,t,x,h}_T=\nabla\phi_n(X_T^{t,x})e^{(T-t)A}Bh,
\end{array}\right.
\end{equation}
The convergence of $(Y^{n,t,x},Z^{n,t,x})$ to $(Y^{t,x},Z^{t,x})$ stated in the proof of the previous
Theorem \ref{teoKolmo} still holds true, as well as estimates (\ref{stimaFnVn}).
In order to investigate the convergence of $F^{n,t,x,h}$, we get
\begin{align*}
  F^{n,t,x,h}_t-F^{j,t,x,h}_t& = \E\left[\nabla\phi_n(X_T^{t,x})-
\nabla\phi_j(X_T^{t,x})\right]e^{(T-t)A}Bh \\
&+ \E\int_t^T \left( 
\nabla\psi_n(Z^{n,t,x}_\tau)V^{n,t,x,h}_\tau-
\nabla\psi_j(Z^{j,t,x}_\tau)V^{j,t,x,h}_\tau\right)\,d\tau.
\end{align*}
and consequently
\begin{align*}
 \vert F^{n,t,x,h}_t-F^{j,t,x,h}_t\vert &\leq \vert\E\left[\nabla\phi_n(X_T^{t,x})-
\nabla\phi_j(X_T^{t,x})\right]e^{(T-t)A}Bh\vert\\
& +\vert\E \int_t^T \left( 
\nabla\psi_n(Z^{n,t,x}_\tau)V^{n,t,x,h}_\tau- \nabla\psi_n(Z^{j,t,x}_\tau)V^{j,t,x,h}_\tau\right)\,d\tau
\vert \\
&+\vert\E \int_t^T \left( 
\nabla\psi_n(Z^{j,t,x}_\tau)V^{j,t,x,h}_\tau- \nabla\psi_j(Z^{j,t,x}_\tau)V^{j,t,x,h}_\tau\right)\,d\tau
\vert=I+II+III
\end{align*}
The terms $I$ and $II$ can be treated as in the proof of Theorem \ref{teoKolmo}, since, by the lipschitz
character of $\psi$ every $\psi_n$ is Lipschitz continuous with the same Lipschitz constant as $\psi$, and the derivative $\nabla \psi_n$ is bounded, uniformly with respect to $n$.
It remains to estimate $III$: similarly to what we have done in the proof of Theorem \ref{teoKolmo}, 
we fix $0<\delta<T-t$ which will be chosen in the following and in order to estimate $III$ we start by estimating
\begin{align*}
& \vert\E \int_{t+\delta}^T \left( 
\nabla\psi_n(Z^{j,t,x}_\tau)V^{j,t,x,h}_\tau- \nabla\psi_j(Z^{j,t,x}_\tau)V^{j,t,x,h}_\tau\right)\,d\tau
\vert \\
&=\vert \int_{t+\delta}^T\int_{H}\left(\psi_n\left(v_j(s,  y+e^{(s-t)A}x)\right)  -
\psi_j\left(v_j(s,  y+e^{(s-t)A}x)\right) \right) \\
& \left\langle Q_{s-t}^{-1/2}%
e^{(s-t)A}Bh,Q_{s-t}^{-1/2}y\right\rangle \mathcal{N}\left(  0,Q_{s-t}\right)  \left(
dy\right)\,ds\vert.\\
\end{align*}
By the dominated convergence Theorem and by the uniform convergence of the sequence $(\psi_n)_n$
this integral goes to $0$ for $n,j\rightarrow \infty$.
We can conclude that, choosen $\varepsilon>0$
\begin{equation}\label{primastimaIII}
 \vert\E \int_{t+\delta}^T \left( 
\nabla\psi(Z^{n,t,x}_\tau)V^{n,t,x,h}_\tau- \nabla\psi(Z^{j,t,x}_\tau)V^{j,t,x,h}_\tau\right)\,d\tau
\vert \leq \varepsilon,
\end{equation}
independently on the choice of $\delta$. We now estimate
\[
\vert\E \int_t^{t+\delta} \left( 
\nabla\psi_n(Z^{j,t,x}_\tau)V^{j,t,x,h}_\tau- \nabla\psi_j(Z^{j,t,x}_\tau)V^{j,t,x,h}_\tau\right)\,d\tau
\vert.
\]
Since $\nabla\psi_n$ and $\nabla\psi_j$ are uniformly
bounded, and $V^{j,t,x,h}$ and $V^{n,t,x,h}$ are uniformly
bounded in $L^2(\Omega \times [0,T])$, so for any
$ \varepsilon >0$ we can choose $\delta$ such that
\[
 \vert\E \int_t^{t+\delta} \left( 
\nabla\psi(Z^{n,t,x}_\tau)V^{n,t,x,h}_\tau- \nabla\psi(Z^{j,t,x}_\tau)V^{j,t,x,h}_\tau\right)\,d\tau
\vert\leq \varepsilon.
\]
and so we can conclude that for $n,j\geq \bar n$
\[
 III\leq2\varepsilon.
\]
We have proved that
\begin{align*}
 \vert F^{n,t,x,h}_t-F^{j,t,x,h}_t\vert &\leq c(T-t)\Vert\phi_n-\phi_j)\Vert\vert h\vert+\varepsilon.
\end{align*} Now the proof goes on as the proof of Theorem \ref{teoKolmo}.
\qed

As a byproduct of the previous proof we have the following identification of $Z^{t,x}_t$ with
$\nabla^B v(t,x)$. 
\begin{corollary}\label{cor:identifZ}
 Let $(Y^{t,x},Z^{t,x})$ be the solution of the BSDE in the FBSDE
(\ref{fbsde}), assume that hypotheses \ref{ip su AB} and \ref{ip costo} hold true and let $v$ be mild solution of
the semilinear Kolmogorov equation (\ref{Kolmo}).
Then 
\begin{equation}\label{identifZ}
 \nabla^B v(\tau,X_\tau^{t,x})= Z_\tau^{\tau,X_\tau^{t,x}},  \quad\P-\text{a.s.}
\text{ for almost all }t\in[0,T].
\end{equation}
\end{corollary}
\dim
Let $v_n$ be the solution of the approximating Kolmogorov equation (\ref{Kolmo_n}) and let
$(Y^{n,t,x},Z^{n,t,x})$ be the solution of the BSDE in the FBSDE
(\ref{fbsde_n}).
We know that for all $t\in[0,T],\,x\in H$ the following identification holds true:
\begin{equation*}
 \nabla^B v_n(t,x)= Z_t^{n,t,x},  \text{ for almost all }t\in[0,T],\,x\in H.
\end{equation*}
Since $\nabla^B v_n(t,x)\rightarrow \nabla^B v(t,x)$ for all $t\in [0,T)$ and $x \in H$,
and $\nabla^B v_n(t,x)$ is bounded uniformly with respect to $n$, we get that
$\nabla^B v_n(\tau,X_\tau^{t,x})\rightarrow \nabla^B v(\tau,X_\tau^{t,x})$ in $L^2(\Omega\times[0,T])$, and so we get the identification (\ref{identifZ}).
\qed

\begin{remark}\label{remark:OUpert-inKolmo} We show how to handle the case of a semilinear Kolmogorov equation (\ref{Kolmo}) where $\call$ is the generator of a perturbed Ornstein-Uhlenbeck process as (\ref{ornstein-pert}),
that is 
\[
(\call_t f)(x)=\frac{1}{2}(Tr BB^* \nabla^2 f)(x)+\<Ax+B F(t,x),\nabla f(x)\>.
\]
We have already noticed in remark \ref{remark:OUpertReg} that if $A$ and $B$ satisfy hypothesis \ref{ip su AB}
and if $F$ is sufficiently regular, then the perturbet Ornstein-Uhlenbeck transition semigroup satisfies the regularizing property stated in Lemma \ref{ipH su fi-cons}.

\noindent We do not use this property for the solution of the semilinear Kolmogorov equation, but an equivalent representation
of the mild solution in terms of an Ornstein-Uhlenbeck transition semigroup based on the Girsanov transform.
To this aim, notice that, at least in the case of $\phi$ and $\psi$ differentiable,
we can apply the Girsanov theorem in the forward-backward system 
\begin{equation*}
\left\{
\begin{array}
[c]{l}%
dX_\tau  =AX_\tau d\tau+BF(\tau,X_\tau) d\tau+BdW_\tau
,\text{ \ \ \ }\tau\in\left[  t,T\right], \\
X_\tau =x,\text{ \ \ \ }\tau\in\left[  0,t\right], \\
\dis
 dY_\tau^{t,x}=-\psi(Z_\tau^{t,x})-l(\tau, X^{t,x}_\tau)\;d\tau+Z^{t,x}_\tau\;dW_\tau,
 \qquad \tau\in [0,T],
  \\\dis
  Y_T^{t,x}=\phi(X_T^{t,x}),
\end{array}
\right.  %
\end{equation*}
or we can follow \cite{Go1}.
We get that the mild solution of equation (\ref{Kolmo}) can be represented, for all $t \in [0,T]$, $x \in H$, as
\begin{align*}
v(t,x)&=R_{t,T}\left[  \phi\right]  \left(  x\right)  +\int_{t}^{T}R_{t,s}\left[l\left(s,\cdot\right)\right(x)ds \\
&+\int_{t}^{T}%
R_{t,s}\left[  \psi\left(\nabla^B 
v\left(  s,\cdot\right)\right)  \right] (  x)  ds
+\int_{t}^{T}R_{t,s}\left[ \nabla^B v(s,\cdot) F (s,\cdot)\right] (  x)  ds. %\label{solmildkolmoeq}%
\end{align*}
Here $(R_{t,T})_{t\in[0,T]}$ is the transition semigroup of the corresponding
Ornstein-Uhlenbeck process
\begin{equation*}
\left\{
\begin{array}
[c]{l}%
dX_\tau  =AX_\tau d\tau+BdW_\tau
,\text{ \ \ \ }\tau\in\left[  t,T\right], \\
X_t =x,\text{ \ \ \ }\tau\in\left[  0,t\right].
\end{array}
\right.
\end{equation*}
The new Hamiltonian function is given by 
\begin{equation}\label{newham}
  \tilde \psi(t,x,z):=\psi(z)+zF(t,x)
\end{equation}
and by a straightforward generalization of Theorem \ref{teoKolmo:gen} the case of a perturbed Ornstein-Uhlenbeck process is covered.
\qed
\end{remark}

\section{The semilinear Kolmogorov equation in the Banach space framework}
\label{sezioneKolmoBanach}
In this Section we revisit the results contained in the previous sections in the case when the problem is considered in a Banach space $E$ instead of the Hilbert space $H$. For the Banach space framework, we mainly refer to the papers \cite{Mas1}, \cite{Mas-Ban} and \cite{Mas-SAP}. From now we assume that $E$ admits a countable Schauder basis.

\noindent From now on, we consider a Banach space $E$ continuously and densely embedded in a real and separable Hilbert space $H$. We have to adequate the Hilbert space case to this new context.
In general, we
can not guarantee that $B\left(  \Xi\right)  \subset E.$ We make the following
assumption which is verified in most of the applications.

\begin{hypothesis}
\label{ip su CSI}There exists a subspace $\Xi_{0}$ dense in $\Xi$ such that
$B\left(  \Xi_{0}\right)  \subset E$
\end{hypothesis}

% \begin{hypothesis}
% \label{ip su CSIbis}There exists  $0<\alpha<\dfrac{1}{2}$ such that for all 
% \[
%  \Vert e^{tA}B\xi\Vert_\E\leq \dfrac{C}{t^\alpha}
% \]
% \end{hypothesis}

In the following we extend the definition of $B$ differentiability to continuous functions $f:E \rightarrow \mathbb{R} $, as in \cite{Mas1}, definition 2.7.

\begin{definition}
\label{def G der banach}For a map $f:E\rightarrow\mathbb{R}$ the
$B$-directional derivative $\nabla^{B}f$ at a point $x\in E$ in direction$\ \xi
\in\Xi_{0}$ is defined as usual.
% \[
% \nabla^{B}f\left(  x;\xi\right)  =\lim_{s\rightarrow0}\frac{f\left(
% x+sG\xi\right)  -f\left(  x\right)  }{s},\text{ }s\in\mathbb{R}\text{.}%
% \]
We say that a continuous function $f$ is $B$-Gateaux differentiable at a point
$x\in E$ if $f$ admits the $B$-directional derivative $\nabla^{B}f\left(
x;\xi\right)  $ in every directions $\xi\in\Xi_{0}$ and there exists a linear
operator $\nabla^{B}f\left(  x\right)  $ from $\Xi_{0}$ with values in
$\R$, such that $\nabla^{B}f\left(  x;\xi\right)  =\nabla^{B}f\left(
x\right)  \xi$ and $\left|  \nabla^{B}f\left(  x\right)  \xi\right|  \leq
C_{x}\left\|  \xi\right\|  _{\Xi}$, where $C_{x}$ does not depend on $\xi$. So
the operator $\nabla^{B}f\left(  x\right)  $\ can be extended to the whole
$\Xi$, and we denote this extension again by $\nabla^{B}f\left(  x\right)  $,
the $B$-gradient of $f$ at $x$. The definition of $f$ $B$-G\^ateaux differential
and of the class $\calg^{B}\left(  E\right)  $ of
functions $f\in C_b(E)$ that are $B$-G\^ateaux differentiable follows in the usual natural way, see \cite{Mas1}.
\end{definition}
%%%%%%%%%%%%%%%%%%
We consider an Ornstein-Uhlenbeck process in $E$, that is a Markov process $X$ solution to equation%
\begin{equation}
\left\{
\begin{array}
[c]{l}%
dX_\tau  =AX_\tau d\tau+BdW_\tau
,\text{ \ \ \ }\tau\in\left[  t,T\right] \\
X_t =x,
\end{array}
\right.  \label{ornsteinBanach}%
\end{equation}
On the operator $A$ and $B$ we make the following assumptions:
\begin{hypothesis}
\label{ip su AB-Banach} We assume that $A$ generates a $C_{0}$ semigroup in $E$;
and we suppose that there exists $\omega
\in\R$ such that $\left\|  e^{tA}\right\|  _{L\left(  E,E\right)
}\leq e^{\omega t}$, for all $0\leq t\leq T$. We assume that $e^{tA}$,
$t\geq0$, admits an extension to a $C_{0}$ semigroup of bounded linear
operators in $H$, whose generator is denoted by $A_{0}$ or by $A$ if no
confusion is possible.
The operator $B\in L\left(  \Xi,H\right)  $ is such that the operators 
\[
Q_{\tau}=\int_{0}^{\tau}e^{sA_{0}}BB^{\ast}e^{sA_{0}^{\ast}}ds\text{, }%
\tau\geq0,
\] are
of trace class for every $\tau\in\left[  0,T\right]  $, so that the stochastic convolution
$W_{A}\left(  \tau\right)  $ is well defined in $H$. We assume further that $W_{A}\left(  \tau\right)$
admits an $E$-continuous version.
\end{hypothesis}
We make the following assumptions on $A$ and $B$, which imply $B$-regularizing properties on the 
Ornstein-Uhlenbeck transition semigroup $P_t$.
\begin{hypothesis}
\label{ipH su fi-Banach}The operators $A$ and $B$ are such that%
\begin{equation}
e^{tA_{0}}B\left(  \Xi\right)  \subset Q_{t}^{1/2}\left(  H\right)
.\label{ornstein inclusione Banach}%
\end{equation}
and
\begin{equation}\label{ornstein inclusione:stima-Banach}
\left\|  Q_{t}^{-1/2}e^{tA_0}B\right\|  \leq c(t)\text{, for }0<t\leq T.
\end{equation}
with $c(t)$ as in Hypothesis \ref{ipH su fi}.
\end{hypothesis}
Under hypothesis \ref{ipH su fi-Banach} it is possible to prove that the transition
semigroup $P_t$ enjoys the $B$-regularizing property, with an anlogous of Lemma \ref{ipH su fi-cons}, following also lemma 2.8 in \cite{Mas1}.
Let us consider also a perturbed Ornstein-Uhlenbeck process in $E$:
\begin{equation}
\left\{
\begin{array}
[c]{l}%
dX_\tau  =AX_\tau d\tau+BF(X_\tau)d\tau+BdW_\tau
,\text{ \ \ \ }\tau\in\left[  t,T\right] \\
X_t =x.
\end{array}
\right.  \label{pert-ornsteinBanach}%
\end{equation}
In order to get the $B$-regularizing property for the transition semigroup of the perturbed Ornstein-Uhlenbeck process in $E$, we have further to assume:
\begin{hypothesis}
\label{ip su CSIbis}For every $\xi\in\Xi$ and for every $t>0$ $\left\|
e^{tA_{0}}B\xi\right\|  _{E}\leq ct^{-\alpha}\left\|  \xi\right\|  _{\Xi}$,
for some constant $c>0$ and $0<\alpha<\dfrac{1}{2}$.
\end{hypothesis}
Then, with the nonlinear term $F$ in (\ref{pert-ornsteinBanach}) lipschitz continuous and differentiable,
following \cite{Mas1}, Section 4, the $B$ regularizing property for the transition semigroup of the perturbed
Ornstein-Uhlenbeck process can be proved, with blow up given by the function $c(t)$ in hypothesis
\ref{ipH su fi-Banach}.

We now solve semilinear Kolmogorov equations in Banach spaces.
Let $\call_t $ be the generator of the transition Ornstein-Uhlenbeck transition
semigroup $P_t$ related to the $E$-valued process defined by (\ref{ornsteinBanach}), that is, at least
formally,
$$
(\call f)(x)=\frac{1}{2}(Tr BB^* \nabla^2 f)(x)+\<Ax,\nabla f(x)\>.
$$
Let us consider the following equation
\begin{equation}
\left\{
\begin{array}
[c]{l}%
-\frac{\partial v}{\partial t}(t,x)=\mathcal{L} v\left(  t,x\right)
+\psi\left(  \nabla^{B}v\left(  t,x\right)  \right)+l(t,x)  ,\text{ \ \ \ \ }t\in\left[  0,T\right]
,\text{ }x\in E\\
v(T,x)=\phi\left(  x\right)  ,
\end{array}
\right.  \label{KolmoBanach}%
\end{equation}
We can adequate immediately definition \ref{solmildkolmo} in order to get the notion of mild
solution of a semilinear Kolmogorov equation in $E$. Also hypothesis
\ref{ip costo} can be immediately adequated to the Banach space framework by substituting
the real and separable Hilbert space $H$ with the Banach space $E$.

Moreover, the counterpart of Theorems \ref{teoKolmo} and \ref{teoKolmo:gen} holds true, and it is summarized in the following Theorem.
\begin{theorem}\label{teoKolmo:genBanach}
 Assume that hypotheses \ref{ip su CSI}, \ref{ip su AB} and \ref{ip costo}, with the Banach space $E$ in the place of the Hilbert space $H$,
hold true.
Then equation (\ref{Kolmo}) admits a unique mild solution $v(t,x)$ according to definition
(\ref{defsolmildkolmo}), and $v(t,x)$ is given by
\[
 v(t,x):=Y^{t,x}_t
\]
with $Y^{t,x}$ solution to the BSDE in (\ref{fbsde}).
\end{theorem}
\dim We give the sketch of the proof, in the case $l=0$ for the sake of simplicity.
We decide not to smooth the final datum $\phi$ by means of the analogous of inf-sup convolutions in Banach spaces, see \cite{Strom}, but to use the approximation procedure introduced in \cite{PZ} for Hilbert spaces and that can be generalized to Banach spaces with a countable Schauder basis. Namely,
we can set
\begin{equation}\label{approxphi}
\phi_{n}\left(  x\right)  =\int_{\mathbb{R}^{n}}\rho_{n}\left(  y-Q_{n}%
x\right)  \phi\left(  \sum_{i=1}^{n}y_{i}e_{i}\right)  dy,
\end{equation}
where $(e_i)_{i\geq 1}$ are the elements of the Schauder basis. We have that
for all $n\geq 1$, $\phi_n$ is lipschitz continuous, with the same lipschitz constant as $\phi$,
and it is G\^ateaux differentiable,
and 
\[
 \vert \phi_n(x)\vert \leq \Vert \phi\Vert_\infty.
\]
We notice that the approximation of $\phi$ by means of $\phi_n$
as defined in (\ref{approxphi}) is only pointwise:
\[
 \lim_{n\rightarrow \infty }\phi_n(x)=\phi(x),\;\forall x \in E.
\]
\noindent Moreover we smooth the Hamiltonian function $\psi$, which is still defined on the Hilbert space $\Xi$ with its inf-sup convolution, see (\ref{infsupconv-psi}). We get an approximating
Kolmogorov equation in the Banach space $E$ given by
\begin{equation}
\left\{
\begin{array}
[c]{l}%
-\frac{\partial v}{\partial t}(t,x)=\mathcal{L} v\left(  t,x\right)
+\psi_n\left(  \nabla^{B}v\left(  t,x\right)  \right)  ,\text{ \ \ \ \ }t\in\left[  0,T\right]
,\text{ }x\in E\\
v(T,x)=\phi_n\left(  x\right)  .
\end{array}
\right.  \label{Kolmo_nBanach}%
\end{equation}
For the approximating Kolmogorov equation (\ref{Kolmo_nBanach}) we can apply \cite{Mas-Ban},
Theorem 6.2, to get that by setting $(Y^{n,t,x},Z^{n,t,x})$ solution of the 
FBSDE (\ref{fbsde_n}) with the process $X$ with values in $E$, the function
$v_n(t,x):=Y^{n,t,x}_t$ turns out to be a mild solution of the approximating Kolmogorov equation (\ref{Kolmo_nBanach}),
for every $n\in \N$. 
We have to show that $v_n(t,x),\nabla^Bv_n(t,x)$ converge and that $v(t,x)=\lim_{n\rightarrow\infty}v_n(t,x)$
is the mild solution of the Kolmogorov equation (\ref{Kolmo}) in $E$.
We let $(Y^n,Z^n)$ be the solution of the approximating BSDE in the FBSDE in the Banach space framework
\begin{equation}\label{fbsde_nBanach}
    \left\{\begin{array}{l}\dis dX_\tau =
AX_\tau d\tau+ B dW_\tau,\quad \tau\in
[t,T]\subset [0,T],
\\\dis
X_t=x,
\\\dis
 dY^n_\tau=-\psi_n(Z^n_\tau)\;d\tau+Z^n_\tau\;dW_\tau,
  \\\dis
  Y^n_T=\phi_n(X_T),
\end{array}\right.
\end{equation}
and $(Y,Z)$ be the solution of the BSDE in the FBSDE in the Banach space framework
\begin{equation}\label{fbsdeBanach}
    \left\{\begin{array}{l}\dis dX_\tau =
AX_\tau d\tau+ B dW_\tau,\quad \tau\in
[t,T]\subset [0,T],
\\\dis
X_t=x,
\\\dis
 dY_\tau=-\psi(Z_\tau)\;d\tau+Z_\tau\;dW_\tau,
  \\\dis
  Y_T=\phi(X_T).
\end{array}\right.
\end{equation}
In order to show that the pair of processes $(Y^n,Z^n)$ converge to the pair $(Y,Z)$, with $\phi_n\rightarrow\phi$ only pointwise, we apply proposition 5.5 in \cite{MasRi}. It is fundamental that $(\phi_n)_n$ are uniformly bounded.
Next, following theorems \ref{teoKolmo} and \ref{teoKolmo:gen}, it is possible to show that
$v_n(t,x),\nabla^Bv_n(t,x)$ converge and that $v(t,x)=\lim_{n\rightarrow\infty}v_n(t,x)$ is the mild solution of the Kolmogorov equation (\ref{Kolmo}) in $E$.
\qed

As in Section \ref{sezioneKolmo} we can notice that $\nabla^B v(t,x)$ is bounded, see Remark \ref{remark:nablaBbounded}, and moreover
as in Corollary \ref{cor:identifZ}, it is possible to show that
\[
\nabla^B v(\tau,X_\tau^{t,x})=Z_\tau^{t,x},\, \P-\text{a.s., for a.a. }t\in[0,T].
\]
Finally, following the outlines given in remark \ref{remark:OUpert-inKolmo}, and assuming also that hypothesis \ref{ip su CSIbis}
holds true, it is possible to handle the case of the second order differential operator $\call_t$ generator of a perturbed Ornstein-Uhlenbeck process in $E$, at least in the case of $F$ lipschitz continuous with respect to $x$.
\section{Solution of the optimal control problem}
\label{applic contr 2}
In this Section we consider the stochastic optimal control problem introduced in Section \ref{applic contr 1}.
This problem can be solved in the Hilbert space framework presented in Section \ref{applic contr 1}, but we skip the solution of the optimal control problem in the Hilbert space framework and we consider it in the Banach space framework, that allows more generality on the choice of the costs, see in particular Section \ref{sez-contr-heat}.

\noindent The controlled state equation takes its values in $E$ and it is given by
\begin{equation}
\left\{
\begin{array}
[c]{l}%
dX^{u}_\tau  =\left[  AX^{u}_\tau +B u_\tau
 \right]  d\tau+BdW_\tau ,\text{ \ \ \ }\tau\in\left[  t,T\right] \\
X^{u}_t  =x.
\end{array}
\right.  \label{sdecontrolforteBanach}%
\end{equation}
The control $u$ is an $\left(\mathcal{F}_{\tau}\right)_{\tau}$-predictable
process with values in a closed and bounded set $U$ of the Hilbert space $\Xi$. Thanks to Hypothesis \ref{ip su AB-Banach} on
the operators $A$ and $B$, and to Hypotheses \ref{ip su CSI} and
\ref{ip su CSIbis}, the mild solution of the controlled state equation (\ref{sdecontrolforteBanach})
is a well defined $E$-valued process.

Beside equation (\ref{sdecontrolforteBanach}), define the cost
\begin{equation}
J\left(  t,x,u\right)  =\mathbb{E}\int_{t}^{T}[
l\left(s,X^{u}_s\right)+g(u_s) ]ds+\mathbb{E}\phi\left(X^{u}_T\right). 
\label{costBanach}%
\end{equation}
for real functions $l$ on $[0,T]\times E$, $\mathbb{\phi}$ on $E$ and $g$ on $U$, which satisfy Hypothesis \ref{ip costo}
with $E$ in the place of $H$.

The control problem will be solved by means of the dynamic programming principle, the value function of the optimal control problem introduced in Section \ref{applic contr 1} is identified with the solution of the HJB equation related. The HJB equation has the same structure of the semilinear Kolmogorov equation we have studied in Section \ref{sezioneKolmoBanach}, equation (\ref{KolmoBanach}), with final datum
$\phi$ equal to the final cost, and with semilinear term $\psi$ equal to the Hamiltonian function defined in (\ref{hamilton}).
\begin{theorem}\label{teo-rel-fond}
 Assume Hypotheses \ref{ip su CSI}, \ref{ip su AB-Banach}, \ref{ipH su fi-Banach}, \ref{ip su CSIbis} and 
\ref{ip costo}, adequated to the Banach space framework, hold true. Let $v$ be the unique mild solution to equation (\ref{KolmoBanach}).
 For every $t\in [0,T]$, $x\in H$ and for all admissible control $u$ we have $J(t,x,u) \geq v(t,x)$, 
 and the equality holds if and only if, for a.a. $s \in [0,T[$, $\mathbb{P}$-a.s.
$$
u_s\in \Gamma\left( \nabla^{B}
v(s ,X^{u,t,x}_s)
\right).
  $$
\end{theorem}
\dim At first we assume that the Hamiltonian function $\psi $ is Lipschitz continuous and G\^ateaux
differentiable and we approximate the current cost $l$, the final cost $\phi$
% and the hamiltoniana function $\psi$
by means of the pointwise approximation introduced in \cite{PZ} and that we have already recalled in the proof of Theorem \ref{teoKolmo:genBanach}, formula (\ref{approxphi}).
 With all the data smooth, we can apply Proposition 5.5, Corollary 5.6 and Theorem 5.7 in \cite{Mas-Ban},
in order to get 
\begin{align}\label{rel-fon_n}
 J_n\left(  t,x,u\right)&:=\mathbb{E}\int_{t}^{T}\left[
l_n\left(s,X^{u}_s\right)+g(u_s) \right]ds+\mathbb{E}\phi_n\left(X^{u}_T\right)\\
 &=v_n(t,x)-\int_t^T\left[\psi\left(\nabla^Bv_n\left(s,X_s\right)\right)-g\left(u_s\right)-
 \nabla^Bv_n\left(s,X_s\right)u_s\right]\,ds,
\end{align}
where $v_n(t,x)$ is the mild solution of the approximating Kolmogorov equation
\begin{equation*}
\left\{
\begin{array}
[c]{l}%
-\frac{\partial v_n}{\partial t}(t,x)=\mathcal{L} v_n\left(  t,x\right)
+\psi\left(  \nabla^{B}v_n\left(  t,x\right)  \right)+l_n(t,x)  ,\text{ \ \ \ \ }t\in\left[  0,T\right]
,\text{ }x\in H\\
v(T,x)=\phi_n\left(  x\right)  ,
\end{array}
\right.  %
\end{equation*}
Letting $n\rightarrow\infty$ in (\ref{rel-fon_n}) we get
\begin{align*}
 J\left(  t,x,u\right)=v(t,x)-\int_t^T\left[\psi\left(\nabla^Bv\left(s,X_s\right)\right)-g\left(u_s\right)-
 \nabla^Bv\left(s,X_s\right)u_s\right]\,ds,
\end{align*}
and by the definition of $\psi$ the result follows.

\noindent If $\psi$ is not G\^ateaux differentiable, we approximate it by means of its inf-sup convolutions $\psi_k$
defined in \ref{infsupconv-psi}.
We get, by Proposition 5.5 in \cite{Mas-Ban} and by Corollary \ref{cor:identifZ}
\begin{align}
\label{relfondk}
 J\left(  t,x,u\right)&=\mathbb{E}\int_{t}^{T}\left[
l\left(s,X^{u}_s\right)+g(u_s) \right]ds+\mathbb{E}\phi\left(X^{u}_T\right)\\
&=\mathbb{E}\int_{t}^{T}\left[
l\left(s,X_s\right)+\psi_k\left(\nabla^Bv_k\left(s,X_s\right)\right)\right]ds+\mathbb{E}\phi\left(X_T\right)\\
&+\E\int_t^T\left[g(u_s)-\psi\left(\nabla^Bv_k\left(s,X_s\right)\right)+\nabla^Bv_k\left(s,X_s\right)u_s\right]\,ds\\
&+\E\int_t^T\left[\psi\left(\nabla^Bv_k\left(s,X_s\right)\right)-\psi_k\left(\nabla^Bv_k\left(s,X_s\right)\right)\right]\,ds\\
 &=v_k(t,x)-\int_t^T\left[\psi\left(\nabla^Bv\left(s,X_s\right)\right)-g\left(u_s\right)-
  \nabla^Bv\left(s,X_s\right)u_s\right]\,ds\\
 &-\int_t^T\left[\psi_k\left(\nabla^Bv_k\left(s,X_s\right)\right)
 -\psi\left(\nabla^Bv_k\left(s,X_s\right)\right)\right]\,ds
\end{align}
where $v_k(t,x)$ is the mild solution of the approximating Kolmogorov equation
\begin{equation*}
\left\{
\begin{array}
[c]{l}%
-\frac{\partial v}{\partial t}(t,x)=\mathcal{L} v\left(  t,x\right)
+\psi_k\left(  \nabla^{B}v\left(  t,x\right)  \right)+l(t,x)  ,\text{ \ \ \ \ }t\in\left[  0,T\right]
,\text{ }x\in H\\
v(T,x)=\phi\left(  x\right)  ,
\end{array}
\right.  %
\end{equation*}
Letting $k\rightarrow\infty$ in (\ref{relfondk}) we get
\begin{align*}\label{rel-fon_n}
 J\left(  t,x,u\right)=v(t,x)-\int_t^T\left[\psi\left(\nabla^Bv\left(s,X^u_s\right)\right)-g\left(u_s\right)-
 \nabla^Bv\left(s,X^u_s\right)u_s\right]\,ds,
\end{align*}
and the Theorem is proved.
\qed

Moreover we can perform the synthesis of the optimal control.
Let us define now the \textit{
optimal feedback law}:
\begin{equation*}%\label{leggecontrolloottima}
u(\tau,x)=\gamma\Big(\nabla^{B}
v(\tau ,X^{u,t,x}_\tau) \Big),\qquad
\tau\in [t,T],\;x\in H.
\end{equation*}
We consider the \textit{closed loop equation} 
\begin{equation}
\left\{
\begin{array}
[c]{l}%
dX^{u}_\tau  =\left[  AX^{u}_\tau +B \gamma(\nabla^{B}
v(\tau ,\overline{X}_\tau)))\,d\tau
 \right]  d\tau+BdW_\tau ,\text{ \ \ \ }\tau\in\left[  t,T\right] \\
X^{u}_t  =x.
\end{array}
\right.  \label{closedloop}%
\end{equation}
and we assume that it
admits a solution
\begin{equation}\label{closedloopmild}
\overline{X}_s= e^{(s-t)A}x_0
+\int_{t}^s e^{(r-t)A}B(\gamma(\nabla^{B}
v(r ,\overline{X}_r)))\,dr
+\int_t^s e^{(r-t)A}B\,dW_r.
\end{equation}
Then the pair $(\overline{u}=u(s,\overline{X}_s),\overline{X}_s)_{s\in[t,T]}$
is optimal for the control problem.
We recall that existence of a solution of the closed loop
equation is not obvious, and that this problem can be avoided by
formulating the optimal control problem in the weak sense, following \cite{FlSo}.
The advantage of the weak formulation is that the closed loop equation is solvable in the weak sense by a Girsanov change of measure.

\subsection{Optimal control problem for the heat equation}
\label{sez-contr-heat}
In this Section we treat stochastic optimal control problems related to stochastic heat equations with control and noise on a subdomain.
Stochastic optimal control problems for reaction diffusion equations
have been extensively studied in the literature.
We cite in particular the papers \cite{Ce1} and \cite{Ce2}, where equations with a more general
structure than equation (\ref{heat equation}) below are treated, but some more smoothing properties on the transition semigroup are required, and the paper \cite{Mas-SAP} where the case of a non linear term in the heat equation is treated and the problem is solved in the Banach space of continuous functions, on the contrary
all the coefficients are asked to be G\^ateaux differentiable.
In the present paper we are able to remove differentiability assumptions on the cost. Moreover we are able to consider only Lipschitz continuous current and final costs in the Banach space of continuous functions, so we are able to treat costs like the supremum of the state.

\noindent Applying the results in section 4, we can also solve the associated HJB equation. We are not able to solve the related HJB equation as in \cite{Mas} by a suitable fixed point argument. Indeed the transition semigroup of the state equation is Strong Feller. This regularizing property is related to the fact that deterministic linear heat equations with control on a subinterval are null controllable, see e.g. \cite{Zu1}, but the minimal energy blows up too fast, see e.g. \cite{FZ} and \cite{Zu}. So since the minimal energy blows up too fast at $0$, namely it is not integrable near $0$, fixed point arguments as the ones used in
\cite{Go1}, \cite{Mas} and \cite{Mas1} cannot be applied to this problem.

In the following we present the controlled heat equation we are able to treat, which is the controlled version of (\ref{heat equation}).

We denote by $H$, which in this case coincides with $\Xi$, the Hilbert space $L^2([0,1])$ and
by $E$ the space of continuous functions on $[0,1]$. We consider the controlled heat equation
\begin{equation}\label{heat equation contr}
 \left\{
  \begin{array}{l}
  \dis
\frac{ \partial y}{\partial s}(s,\xi)= \Delta y(s,\xi)+ 1_{\calo}(\xi)u(s,\xi)
+ 1_{\calo}(\xi)\frac{ \partial W
}{\partial s}(s,\xi), \qquad s\in [t,T],\;
\xi\in [0,1],
\\\dis
y(t,\xi)=x(\xi),
\\\dis
 \dfrac{\partial}{\partial\xi}y(s,\xi)=0, \quad \xi=0,1.
\end{array}
\right.
\end{equation}
Here $\frac{ \partial W
}{\partial s}(s,\xi)$ is a space time white noise and $u:[0,T]\times[0,1]\rightarrow \R$ is the control
process such that $u\in L^2([0,T],L^2([0,1]))$.

\noindent We are able to treat cost functionals
\begin{equation}
J\left(  t,x,u\right)  =\mathbb{E}\int_{t}^{T}\left[l(s,X^u(\cdot))+ g(u_s)\right] ds+\mathbb{E}%
 \phi\left(  X^{u}\left( \cdot\right)  \right)
  , \label{cost-heat}%
\end{equation}
where 
\[
 l:[0,T]\times C([0,1])\rightarrow \R,\quad \phi: C([0,1])\rightarrow \R.
\]
For example, we are able to handle the case of final and current cost given by the supremum of the state over the interval $[0,1]$:
\[
 l(s,x(\cdot))=\sup_{\xi\in[0,1]} x(\xi),\qquad \phi(x(\cdot))=\sup_{\xi\in[0,1]} x(\xi)
\]
To fit our hypotheseis \ref{ip costo} we have to suppose only lipschitz continuity of $\phi$ and of $l$ with respect to $x\in C([0,1])$.

\begin{hypothesis}
\label{heatipotesi}The costs  $l$ and $\phi$ are all Borel
measurable and real valued. Moreover
\begin{enumerate}
\item $ l:\left[  0,T\right]  \times C(\left[  0,1\right])%
\rightarrow\mathbb{R}$ is continuous and bounded, and for all $t\in[0,T]$, 
$l(t,\cdot):C(\left[  0,1\right])%
\rightarrow\R$ is lipschitz continuous.

\item $ \phi:C\left(\left[  0,1\right]\right)  \rightarrow\mathbb{R}$ is
continuous, bounded, and lipschitz continuous.
\end{enumerate}
\end{hypothesis}

The heat equation (\ref{heat equation}%
)\ can be written in abstract way in the Banach space $E$ as
\begin{equation}
\left\{
\begin{array}
[c]{l}%
dX_{\tau}^{u}=  AX_{\tau}^{u}+Bu_{\tau}d\tau+BdW_{\tau}\text{\ \ \ }\tau\in\left[
t,T\right] \\
X_{t}^{u}=x_{0},\\
\end{array}
\right.  \label{heat eq abstract}%
\end{equation}
where $A$ is the Laplace operator in $E=C([0,1])$ with Neumann boundary conditions,
and $B$ is the multiplication by the indicator function $1_{\calo_0}$. The operator $A$
turns out to be the generator of a strongly continuous semigroup in $E$, with an extension to $H$. 
Also hypothesis \ref{ip su CSI} is satisfied by taking $\Xi
_{0}=\left\lbrace f\in C\left(  \left[  0,1\right]  \right):f(a)=f(b)=0 \right\rbrace $, where we recall that $\mathcal O=[a,b]$. With this choice, $\Xi_{0}$ is dense in $\Xi=L^2([0,1])$, and $B (\Xi_0) \subset E=C([0,1])$.
We refer also to \cite{Mas-SAP}
for more details in the reformulation.

We consider the Hamilton Jacobi Bellman equation relative to
(\ref{heat eq abstract})%
\begin{equation}
\left\{
\begin{array}
[c]{l}%
-\frac{\partial v}{\partial t}(t,x)=\call v\left(  t,x\right)
+\psi\left(\nabla^B v\left(  t,x\right)  \right)+l(t,x)  ,\text{ \ \ \ \ \ }%
t\in\left[  0,T\right]  ,\text{ }x\in H,\\
u(T,x)=\phi\left(  x\right)  ,
\end{array}
\right.  \label{heat hamilton}%
\end{equation}
where $\psi$\ is the Hamiltonian function defined in (\ref{hamilton}). In
order to find mild solutions of the
HJB equation and  to solve the optimal control problem, we see that $A$ and $B$ are such that
hypotheses \ref{ip su AB-Banach},  \ref{ip su CSI}, \ref{ip su CSIbis} and \ref{ipH su fi-Banach}
are satisfied. By hypotheses \ref{heatipotesi} on $l$ and $\phi$ in the definition, hypothesis \ref{ip costo} is satisfied.

\begin{theorem}
\label{teo calore}Assume that hypothesis \ref{heatipotesi} holds true. Then equation (\ref{heat hamilton}) has a unique mild solution $v$ and for
all the admissible controls $J\left(
t,x,u\right)  \geq v\left(  t,x\right)  $. Moreover
$J\left(  t,x,u \right)  =v\left(  t,x\right)  $ if
and only if for a.a. $s \in [0,T[$, $\mathbb{P}$-a.s.
$$
u_s\in \Gamma\left( \nabla^{B}
v(s ,X^{u,t,x}_s)
\right).
  $$
\end{theorem}

\end{document}